\documentclass{amsart}
\usepackage{amssymb,latexsym}

\theoremstyle{plain}
\newtheorem{theorem}{Theorem}

\newtheorem{remark}{Remark}

\newtheorem{example}[theorem]{Example}

\theoremstyle{definition}
\newtheorem{definition}{Definition}

\theoremstyle{remark}

\numberwithin{equation}{section}

\begin{document}
\title[Modern Set]{Modern Set}
\author{Jun Tanaka}
\address{University of California, Riverside, USA}
\email{juntanaka@math.ucr.edu, yonigeninnin@gmail.com, junextension@hotmail.com}

\date{May, 25, 2008}

\maketitle


\section{Introduction}

In this paper, we intend to generalize the classical set theory as much as possible. we will do this by freeing sets from the regular properties of classical sets; e.g., the law of excluded middle, the law of non-contradiction, the distributive law, the commutative law,etc.... The fuzzy set theory succeeded in freeing sets from the law of excluded middle and the law of contradiction. However, in order to extend our language, it is more or less unreasonable to keep the commutative law and the distributive law. This modern idea of sets keeps the concept of membership functions but their value are not necessarily in [0,1]; nor do these modern sets form a lattice necessarily. Especially noteworthy is that modern sets are more general than generalized fuzzy set, (Please refer to Nakajima \cite{Naka}). Here is the hierarchy of generality: Modern sets $\geq$ L-fuzzy set $\geq$ Generalized Fuzzy set $\geq$ fuzzy set $\geq$ classical set. These modern sets become the classical sets under the restriction to $O_{x}$, $I_{x}$ for all x in X, and no further conditions are required

The world of natural languages can not be ruled by a single classical logic because the real physical world consists of many aspects, each of which obeys a different logic. As is well-known, the distributive law does not always hold in every logic system. For example, in classical and intuitionistic logic it holds, but in quantum logic it fails.

The Commutative law seems invalid when we talk about anything with temporal order in our daily conversation since every event in time is not reversible. "He is a student and a male". I can also say, "He is a male and a student" but mean the same thing. The two significations commute of course. "I have a health insurance and a car insurance". In this case, the two significations commute. But how about this case "He went to a supermarket and then a drug store". In this case, the two significations do not commute anymore since there is temporal order. Furthermore, we have no expression for "He went to a supermarket and then a drug store" in classical logic or fuzzy logic. In other words, the properties from classical logic, which we take for granted, do not express our thoughts that well.

\section{Preliminaries : Extension of Lattice}

In this section, we shall briefly review the well-known facts about lattice theory (e.g. Birkhoff \cite{Birk}, Iwamura \cite{Iwamura} ), propose an extension lattice, and investigate its properties. (L,$\wedge  $,$\vee  $) is called a lattice, if it is closed under operations $\wedge  $ and $\vee  $, and satisfies, for any elements x,y,z in L:

(L1) the commutative law: x $\wedge$ y = y $\wedge$ x and x $\vee$ y = y $\vee$ x

(L2) the associative law:
\[
\begin{aligned}
x \wedge ( y \wedge z ) = (x \wedge y) \wedge z \ \ \text{and} \ \   x \vee ( y \vee z ) = (x \vee y) \vee z
\end{aligned}
\]

(L3) the absorption law:    x $\vee$ ( y $\wedge$ x ) =x  and   x $\wedge$ ( y $\vee$ x ) = x.

Hereinafter, the lattice (L,$\wedge  $,$\vee  $) will often be written L for simplicity.

A mapping h from a lattice L to another $L'$ is called a lattice-homomorphism, if it satisfies
\[
\begin{aligned}
h ( x \wedge y ) = h(x) \wedge h(y)  \ \ \text{and} \ \   h ( x \vee y ) = h(x) \vee h(y) , \forall x,y \in L.
\end{aligned}
\]

If h is a bijection, that is, if h is one-to-one and onto, it is called a lattice-isomorphism; and in this case, $L'$ is said to be lattice-isomorphic to L.

A lattice (L,$\wedge  $,$\vee  $) is called distributive if, for any x,y,z in L,

(L4) the distributive law holds:
\[
\begin{aligned}
x \vee ( y \wedge z ) = (x \vee y) \wedge (y \vee z)  \ \ \text{and}  \  \   x \wedge ( y \vee z ) = (x \wedge y) \vee ( y \wedge z)
\end{aligned}
\]

A lattice L is called complete if, any subset A of L, L contain the supremum $\vee$ A and the infimum $\wedge$ A. If L is complete, then L itself includes the maximum and minimum elements, which are often denoted by 1 and 0, or I and O, respectively.

\begin{definition} {\textbf{Complete Heyting algebra (cHa)}}

A complete lattice is called a complete Heyting algebra (cHa), if
\[
\begin{aligned}
\vee_{i \in I} \ ( x_{i} \wedge y ) =  (\vee_{i \in I} \  x_{i} ) \wedge y
\end{aligned}
\]
holds for $\forall x_{i} , y \in L $ ($i \in I$); where I is an index set of arbitrary cardinal number.
\end{definition}

A distributive lattice is called a Boolean algebra or a Boolean lattice, if for any element x in L, there exists a unique complement $x^C$ such that x $\vee$ $x^C$ = 1 and $x \wedge x^C$ = 0.

It is well-known that for a set E, the power set P(E) = $2^E$. The set of all subsets of E, is a Boolean algebra.

\section{Preliminaries : Generalized fuzzy sets}
In this section we will consider an algebraic structure of a family of fuzzy sets. We will show that the following three families are mutually equivalent: a ring of generalized fuzzy subsets, an extension of (Boolean) Lattice P(X), and a set of L-fuzzy sets (introduced by Gouen [3]).

\begin{definition}\label{def:1}
A family GF(X), which is closed under operations $\vee$ and $\wedge$, is called a ring of generalized fuzzy subsets of X, if it satisfies:

(1) GF(X) is a complete Heyting algebra with respect to $\vee$ and $\wedge$,

(2) GF(X) contains P(X) = $2^X$ as a sublattice of GF(X),

(3) the operations $\vee$ and $\wedge$ coincide with set operations $\bigcup$ and $\bigcap $, respectively, in P(X), and

(4) for any element A in P(X), A$\vee$X= X and A$\wedge \emptyset = \emptyset$.

\end{definition}

\begin{definition}
Let $L_{x}$ be a lattice which is assigned to each x in X, and let L denote $\{ L_{x} \mid x \in X    \} $. An L-fuzzy set A is characterized by an L-valued membership function $\mu_{A}$ which associates to each point x in X an element $\mu_{A} (x)$ in $L_{x}$.

\end{definition}

\begin{theorem}
If each $L_{x}$ is a cHa, then LF(x), the family of all L-fuzzy sets, is a ring of generalized fuzzy subsets of X.

We will show that LF(X) satisfies all the conditions in Definition $\ref{def:1}$. First, by the definition, LF(X) is closed under $\wedge$ and $\vee$. (1) LF(X) is a cHa, because each $L_{x}$ is a cHa. Condition (2) LF(X) $\supseteqq$ P(X) and (3) $\vee$ = $\cup$ and $\wedge$ = $\cap$ in P(X) follow from the fact that any element of P(X) is defined with the element of LF(X) whose membership function takes just two values, 1 and 0. (4) It follows from $\mu_{X}(x)$=1 and $\mu_{\emptyset}(x)$ = 0 that A $\cup$ X = X and A $\cap $ $\emptyset$ = $\emptyset$, for any A in LF(X).

\end{theorem}

\section{Modern Sets}

\begin{definition} {Boolean Algebra}

In conformity with Birkhoff's book \cite{Birk}, the fundamental operations of intersection and union of elements will be defined by
\[
\begin{aligned}
 x \ \cap \ y &    \ \ \ \ \ \ \text{intersection}& \\
 x \ \cup \ y &    \ \ \ \  \ \ \text{union}&\\
\end{aligned}
\]
As is well known, it follows from the definition of Boolean algebra that there exists a unit element I and a null element O for which we have the following:

\[
\begin{aligned}
 x \cap I = x \ \text{and} & \  x \cap O = O    \\
 x \cup I = I  \ \text{and} & \   x \cup O = x \ \ \ \forall x \in X \\
\end{aligned}
\]
Note that $\cap$ and $\cup$ commute at this point.

\end{definition}

\begin{definition} {\textbf{Weak Boolean algebra }} \label{def:WBA}

Let H be an algebraic space with two distinct operators $\ast_{\wedge}$ , $\ast_{\vee}$ from H to itself. H is called a Weak Boolean Algebra if $\exists $ distinct O, I $\in$ H such that
\[
\begin{aligned}
 O \ast_{\wedge} I = O \ \text{and} & \  I \ast_{\wedge} O = O    \ \ \ \ \ \ \text{O and I commute} \\
 O \ast_{\wedge} O = O \ \text{and} & \   I \ast_{\wedge} I = I   \\
 O \ast_{\vee} I = I \ \text{and} & \  I \ast_{\vee} O = I   \ \ \ \  \ \ \text{O and I commute} \\
 O \ast_{\vee} O = O \ \text{and} &  \ I \ast_{\vee} I = I   \\
\end{aligned}
\]

Please note that $\ast_{\wedge}$ $\ast_{\vee}$ are associated with $\wedge$ and $\vee$, respectively. O and I are associated with the minimum element and the maximum element, respectively.

\end{definition}

\begin{definition} {\textbf{Modern Sets}}

Suppose a weak Boolean Algebra $H_x$ is assigned to each x in X and let H denote $\{ H_x | x  \in X  \}$. Each modern set A is characterized by a membership function $\mu_A$ such that for each x $\in$ X, $\mu_A$ assigns an element $\mu_A (x) \in H_x$. We define $H(X)$ as the family of all modern sets. When A is a set in the ordinary sense of the term (in P(X)), its membership function can take on only two values $O_x$ and $I_x$ with $\mu_A (x)$ = $O_x$ or $I_x$ according to whether x does or does not belong to A. The operations $\ast_{\vee_x}$ and $\ast_{\wedge_x}$ coincide with the set operations $\bigcup$ and $\bigcap $, respectively, in P(X).

A modern set is empty if and only if its membership function is identically $O_x$ for all x in X.

Two modern sets A and B are equal if and only if $\mu_A (x)$ = $\mu_B (x)$ for all x in X.

If $H_x$ is partially ordered by $\leq$ for each x, then we can define containment as follows: A is contained in B if and only if $\mu_A (x) \leq \mu_B (x)$ for all x in X.

$\textbf{Union}$. The union of two modern sets A and B with respective membership functions $\mu_A (x)$ and $\mu_B (x)$ is a modern set, written as C = A $\vee$ B, whose membership function is related to those of A and B by
\[
\mu_C (x) = \mu_A (x) \ast_{\vee_x}  \mu_B (x) , \ x \in X
\]

Note that the order does matter if $H_x$ is not commutative.

$\textbf{Intersection}$ The intersection of two modern sets A and B with respective membership functions $\mu_A (x)$ and $\mu_B (x)$ is a modern set, written as C = A $\wedge$ B, whose membership function is related to those of A and B by
\[
\mu_C (x) = \mu_A (x) \ast_{\wedge_x}  \mu_B (x) , \ x \in X
\]

The notion of complement was not given in Definition \ref{def:WBA} since we can trivially define the complement on a Modern Set for each x in X such as $ \cdot {}^C : H_x \rightarrow H_x  $ where
\[ (\{ O_x \}) ^C = I_x \ \ \text{and} \ \ (\{ I_x \}) ^C = O_x \] even while $(\{ A_x \}) ^C$ can be anything for all $A_x \in H_x , A_x \neq O_x, I_x$ such that $((\{ A_x \}) ^C ) ^C$ =  $\{ A_x \}$

\end{definition}

\begin{theorem}
For any modern sets A, B , C whose membership function take on values $O_x$ or $I_x$, (L1)-(L4) hold.

\begin{proof}
Since the proof is more or less clear, herein we will briefly indicate the proof for the commutative and distributive cases. Let's check the commutative law. Call C = A $\vee$ B.   We only have to check the four possible cases:
\[
\begin{aligned}
 O_x \ast_{\wedge} I_x = O_x  \ \text{and} & \  I_x \ast_{\wedge} O_x = O_x     \\
 O_x \ast_{\wedge} O_x = O_x \ \text{and} & \   I_x \ast_{\wedge} I_x = I_x  \ \ \text{for all x} \ \in X \\
\end{aligned}
\]
Thus, $ \mu_C (x) = \mu_A (x) \ast_{\vee_x}  \mu_B (x) $ =  $ \mu_B (x) \ast_{\vee_x}     \mu_A (x) $. Therefore, C = A $\vee$ B = B $\vee$ A.

We can show the commutative law for intersection similarly.

We will show the two most important cases of the distributive law briefly. Call C = A $\vee$ ( B $\wedge$ C ).
\[
\begin{aligned}
&O_x \ast_{\vee} ( I_x \ast_{\wedge} O_x  )= O_x  \ast_{\vee} O_x  = O_x  = I_x  \ast_{\wedge} O_x = (  I_x  \ast_{\vee} O_x  )  \ast_{\wedge} (  O_x  \ast_{\vee} O_x    )  \\
&O_x \ast_{\vee} ( I_x \ast_{\wedge} I_x  )= O_x  \ast_{\vee} I_x  = I_x  = I_x  \ast_{\wedge} I_x = (  O_x \ast_{\vee} I_x  )  \ast_{\wedge} (  O_x  \ast_{\vee} I_x    ) \\
 \text{for all x} \ \in X \\
\end{aligned}
\]
Please check the rest of cases. That should not be too difficult.

\end{proof}
\end{theorem}

\begin{theorem}
$H_x$ is commutative for all x in X iff H = $\{H_x  \ |  x  \in X   \}$ is commutative.
\begin{proof}
It is obvious.
\end{proof}
\end{theorem}

\begin{theorem}
$H_x$ is distributive for all x in X iff H = $\{H_x  \ |  x  \in X   \}$ is distributive.
\begin{proof}
It is obvious.
\end{proof}
\end{theorem}

\begin{theorem}
Similar equivalent statement (as above) hold for the absorption law, the law of excluded middle, the law of non-contradiction, the associative law, etc...
\end{theorem}

\begin{remark}
As you see from the previous theorems, if ($H_x$, $\ast_{\vee_x}$, $\ast_{\wedge_x}$, $\cdot {}^C $ ) is defined to be in the sense of fuzzy sets such as $\ast_{\vee_x}$ = max, $\ast_{\wedge_x}$ = min, $( \cdot )^C $ = 1- ($\cdot $) in the interval [0,1], then $H(X)$ becomes the family of fuzzy sets.

If $H_x$ satisfy the definition of Lattice for each x in X , then $H(X)$ becomes a family of L-fuzzy set. If $H_x$ additionally satisfy the definition of a complete Heyting algebra, it is a ring of generalized fuzzy subsets.
\begin{proof}
The proof is clear.
\end{proof}
\end{remark}

\begin{example}\label{Ex:1}
We will present an example of a non-commutative modern set. Let $H_x$ be a space of linear bounded operators on a Hilbert space $    \mathcal{H}$ for each x in X. Then we take the zero $O_x $ and the identity $I_x$ in $H_x$. We define the composition $\circ$ =  $\ast_{\wedge_x}$ and the addition + = $\ast_{\vee_x}$. In order to create a weak Boolean Algebra, we must define an equivalence class $\sim$ as $A \sim B$ iff A = B or there exist n, m $\in \mathbb{N}-\{ 0 \}$ such that A = $nI_x$ and B = $mI_x$.

\[
\begin{aligned}
 O_x ( I_x ) = O_x \ \text{and} & \  I_x ( O_x ) = O_x    \ \ \ \ \ \  O_x \ \text{and} \ I_x \ \text{commute} \\
 O_x ( O_x ) = O_x \ \text{and} & \   I_x ( I_x ) = I_x   \\
 O_x + I_x = I_x \ \text{and} & \  I_x + O_x = I_x   \ \ \ \  \ \ O_x \ \text{and} \ I_x \ \text{commute} \\
 O_x + O_x = O_x \ \text{and} &  \ I_x + I_x = 2I_x=I_x   \\
\end{aligned}
\]

 Now we have a weak Boolean Algebra $H_x$ where the composition operation is typically not commutative. Thus H(X) is a family of modern sets of non-commutative type if one of $H_x$ is not commutative under the composition. Needless to say, if all of $H_x$ is commutative under the composition, then H(X) satisfies the commutative law.
\end{example}

\begin{example}
Let (M, n$\times$n , + , $\cdot$) be an n$\times$n matrix space closed under addition + and matrix multiplication $\cdot$. Then call the matrix multiplication identity I and the addition identity matrix O. Now we take $\cdot$ = $\ast_{\wedge_x}$ and + = $\ast_{\vee_x}$. By considering the same equivalence class as in Example \ref{Ex:1}, we can create a modern set from the n$\times$n matrix space.
\end{example}

\begin{theorem}{Gelfand Theorem}\label{The:GT}

If $\mathfrak{U}$ is a commutative a $C^{\ast}$-algebra, then $\mathfrak{U}$ is $\ast$-isomorphism to $\mathbf{C}(X)$, some compact Hausdorff space X.
\end{theorem}

\begin{remark}
For a given commutative $C^{\ast}$-algebra, we have a representation of it as a space of continuous functions on some compact Hausdorff by theorem \ref{The:GT}. Thus we can always create a commutative Modern set from any commutative $C^{\ast}$-algebra under the same construction as in Example \ref{Ex:1}.
\end{remark}

We give the following definition and theorem as more examples of Modern Sets.

\begin{definition}
By a representation of a $C^{\ast}$-algebra $\mathfrak{U}$ on a Hilbert space $    \mathcal{H}$, we mean a  $^{\ast} $ homomorphism $\varphi$ from $\mathfrak{U}$ into $ \mathfrak{B} (   \mathcal{H} )$. If in addition, $\varphi$ is one-to-one, it is called a faithful representation.
\end{definition}

\begin{theorem}{\textbf{The Gelfand-Neumark Theorem}}\label{The:GNT}

Each $C^{\ast}$-algebra has a faithful representation on some Hilbert space.

\end{theorem}

\begin{example}
For a given $C^{\ast}$-algebra, we have a representation of it as bounded linear operators on a Hilbert space by theorem \ref{The:GNT}. Thus we can always create a Modern set from any $C^{\ast}$-algebra under the same construction as in Example \ref{Ex:1}.
\end{example}

\section{Conclusion and observation}
As we mentioned in the Introduction, systematic expression of our thought requires room for at least the non-commutative property. I strongly believe that this new logic system will open up a new blanch of Artificial Intelligence. Property-like verbs such as "be", "have", and "own" seem valid in classical logic. However, most of the other verbs are required to be non-commutative with respect to objects and time. This modern set does not need to be commutative, in some sense, this is closer to the system of our thought. We need further investigation to improve the systematic expression of our thought in order to create a real Artificial Intelligence. I dream the day will come, when we make a real AI.




\begin{thebibliography}{99}



\bibitem{Birk}
G. Birkhoff, Lattice Theory, 3rd ed. AMS colloquim Publication, Providence, RI, 1967.

\bibitem{Lfuzzy}
J.A. Goguen, L-fuzzy sets, J. Math. Anal. Appl., 18, 145-174, 1967

\bibitem{being and time}
M. Heidegger, being and time, HarperOne; Revised edition, 1962 (English)

\bibitem{Iwamura}
T. Iwamura, Sokuron (Kyoritsu Shuppan, Tokyo, 1966).


\bibitem{operator algebra}
R. V. Kadison, J.R. Ringrose, Fundamentals of the Theory of Opreator Algebra, AMS, 1997.

\bibitem{Naka}
N. Nakamura, generalized fuzzy sets, Fuzzy sets and Systems 32 (1989), 307-314.

\bibitem{Take}
G. Takeuti, Senkei-Daisu to ryoushi-rikigaku, 131-162, Shokabo, Tokyo, 1981 (Japanese)


\bibitem{Zad}
L.A. Zadeh, Fuzzy sets, Information and Control 8 (1965), 338-353.


\end{thebibliography}
\end{document}